\documentclass{article}
\usepackage{graphicx} 
\usepackage{xcolor}

\usepackage{amsmath,amsfonts,amssymb,amsthm}
\title{A note on dispersing particles on a line}
\author{Alan Frieze\thanks{Research supported in part by NSF grant DMS1362785
} and Wesley Pegden\thanks{Research supported in part by NSF grant DMS1363136}\\Department of Mathematical Sciences\\Carnegie Mellon University\\Pittsburgh PA 15213}

\usepackage[shortlabels]{enumitem}

\newtheorem{theorem}{Theorem}

\newtheorem{lemma}[theorem]{Lemma}

\newtheorem{observation}[theorem]{Observation}

\theoremstyle{definition}
\newtheorem*{remark}{Remark}

\newcommand{\1}{\mathbf{1}}

\newcommand{\Z}{\mathbb{Z}}
\newcommand{\visits}{\mathrm{arrivals}}
\newcommand{\fvisits}{\mathrm{firstarrivals}}
\newcommand{\svisits}{\mathrm{soloarrivals}}
\newcommand{\E}{\mathbb{E}}
\begin{document}

\maketitle

\begin{abstract}
We consider a synchronous dispersion process introduced in \cite{CRRS} and we show that on the infinite line the final set of occupied sites takes up $O(n)$ space, where $n$ is the number of particles involved.
\end{abstract}

\section{Introduction}
This note concerns a synchronous dispersion process introduced by Cooper, McDowell, Radzik, Rivera and Shiraga \cite{CRRS}. In their model, configurations of particles on the vertices of a graph evolve in discrete time; at each time step, each particle at a vertex with at least 2 particles in total moves to a neighbour chosen independently and uniformly at random; we say the vertex \emph{topples}. (The dispersion process thus ends at the first step when each vertex has at most one particle.) In \cite{CRRS}, they study the behavior of this process on various graphs when begun from a configuration consisting of $n$ particles at one vertex of the graph, with all other vertices initially empty.

They studied this process on a variety of graphs, one of which was the two-way infinite path $L$, with vertex set $\Z$ and edge set $\{\{i,i+1\}, i\in \Z\}$.  They proved that if the initial configuration of particles consists just of $n$ particles at the origin 0, then w.h.p. the furthest particle from the origin is at distance $O(n\log n)$ when the process stops. In this note we reduce this to $O(n)$.
\begin{theorem}\label{disperse}
Suppose that we begin the dispersion process on $L$ with $n$ particles at the origin. Then there is an absolute constant $c>0$ such that w.h.p.~the furthest particle from the origin is at distance at most $cn$ when the process stops.
\end{theorem}
The proof we give here does not depend on synchronous topples; in fact, we can allow an adversary to choose, at each time $t$, a subset of vertices all with at least two particles, and then topple all vertices in that subset.

\begin{remark}
    This version corrects \cite{olddisperse}.  In particular, given a biased random walk to the left, letting an adversary choose when the walk takes a step cannot increase the probability that the walk makes it to some $x>0$, nor can it increase the expected number of arrivals to $x>0$, but it is not true that it cannot increase the probability that the walk is at $x>0$ at some fixed time, a hypothesis that was implicitly used in the claim of equation (3) in \cite{olddisperse}.  The proof here starts from the same perspective---analyzing the same ordered process as in \cite{olddisperse}, defined below, and using the fact that gaps between particles in this process tend to contract---but is written to bound and leverage the expected number of ``arrivals'' gaps make to a large value, rather than the probability that a gap is at a large value at a particular time.  We give a complete self-contained proof here, so that no reference to \cite{olddisperse} is required.
\end{remark}

\section{Proof of Theorem \ref{disperse}}

For our proof, we will study another equivalent process; which we call the \emph{ordered} disperson process on $L$.  In the ordered dispersion process, each particle $p_j$ has an assigned label $j$ from $1$ to $n$. We let $p_{j,t}$ denote the position of particle $j$ at time $t$, and $N_{x,t}$ denote the number of particles at position $x$ at time $t$.  To compute one time step of the ordered process, the rules for the original dispersion process can be applied, and then the particles simply relabeled so that $j_1<j_2$ implies $p_{j_1,t}\leq p_{j_2,t}$.

Of course, it is also possible to characterize the ordered process in more complicated way without relabeling. In this view, the probability that a particle moves left or right when its vertex topples is rarely $\tfrac 1 2$, and, in general, depends both on the number of other particles in its stack, as well as the number of particles occupying the vertices adjacent to its vertex, and whether those vertices topple.  Our proof works by analyzing this more direct (and more complicated) view of the ordered process.

A key advantage to this ordered process is that in the ordered process, particles are biased to move towards empty gaps.  In particular, we have the following:

\begin{observation}
    \label{bias}
    Suppose at time $t-1$ we have that $p_{i,t-1}\geq p_{i-1,t-1}+2$.  If $p_{i,t-1}$ topples at time $t-1$, then we have that
    \[
    p_{i,t}=p_{i,t-1}-1\quad\text{with probability at least }\frac 3 4.
    \]
\end{observation}
\begin{proof}
 We have $\Pr(p_{i,t}=p_{i,t-1}-1)=1-2^{-a}$ where $a=N_{p_{i,t-1},t-1}$.
\end{proof}
  The basic idea behind our proof is to show, by induction on $j$, that particle $p_j$ arrives at sites $x$ far to the right very few times in expectation.  (It is important to note that we are only bounding the number of \emph{arrivals} at $x$, not the total time spent there.)  In particular, we can imagine watching $p_j$ after each time it separates from $p_{j-1}$ at some vertex $x_{j-1}$.  So long as it remains distance 2 from $p_{j-1}$, Observation \ref{bias} ensures that its steps behave like a left-biased random walk.  On the other hand, the following observation ensures that once it is within distance 1 of $p_{j-1}$, it has a reasonable chance of visiting it.  Together these two facts will imply that each time $p_j$ and $p_{j-1}$ separate at a vertex $x_{j-1}$, $p_j$ visits sites $x\gg x_j$ few times in expectation before its next encounter with the particle $p_{j-1}$, as will be encapsulated in Lemma \ref{walkvisits}, below.
\begin{observation}
    \label{ends}
        Suppose at time $t-1$ we have that $p_{i,t-1}=p_{i-1,t-1}+1$.  If $p_{i,t-1}$ topples at time $t-1$ then we have that
    \[
    p_{i,t}=p_{i,t-1}\quad\text{with probability at least }\frac 1 4.
    \]
\end{observation}
\noindent The proof is a straightforward calculation and we defer it to Section \ref{remaining}.

The following lemma is what we use to analyze the walk of a particle $p_j$ from its last encounter with $p_{j-1}$ until its next one.  The adversary plays the role of the aspects of the dispersion process we are not controlling when analyzing the walk of $p_j$; for example, we will not constrain the movement of $p_{j-1}$ (which affects when $p_j$ may revisit $p_{j-1}$) or particles $p_k$, $k>j$ (which affect, for example, with what bias $\geq 3/4$ the particle $p_j$ moves left).
\begin{lemma}
\label{walkvisits}
    Let $p>\frac 1 2$, $r>0$ and consider a process $\{S_n\}$ on $\Z$ which begins with the token at $S_0=0$, and where at every step an adversary can choose to either:
    \begin{enumerate}
    \item Randomly move the token left, right, or keep it in place; $S_{t+1}=S_t+\xi_t$ for $t\in \{-1,0,+1\}$. The probability $\Pr(\xi_t=-1)$ it moves left must be at least $p$. OR: \label{step}
    \item Take an arbitrary step $S_{t+1}=S_t+\xi$ $(\xi_t\in \{-1,0,+1\})$ to the left, right, or in place, but must also choose choose a random $\epsilon_t=\mathrm{Bernoulli}(r)$ and end the walk if $\epsilon_t=1$.
    \label{term}
    \end{enumerate}
    Then with $\rho$ and $c$ as in \eqref{n2} below, no matter how the adversary plays,
    \begin{equation}\label{adversary}
    \forall x\in \Z, \E\Big(\big|\{n\mid S_n=x\}\big|\Big)\leq \rho^{\max(0,x)-c}.
    \end{equation}
\end{lemma}
\noindent Note the bound is of the form one would expect for a left-biased random walk in the absence of the adversary, so the content here is just that the adversary does not have enough power to affect things too much.  The proof is straightforward and we defer it until Section \ref{remaining}. We will apply the lemma in the particular case where 
\begin{equation}\label{n1}
p=\frac 3 4\text{ and } r=\frac 1 4.
\end{equation}
For these values the proof we give shows we can take $\rho=\sqrt{3}/2$, $c=59$.  In particular, when we apply the Lemma the walk $S_n$ will start at $1$, and then we will take
\begin{equation}\label{n2}
\rho=\frac{\sqrt 3}{2}\text{ and }c=60.
\end{equation} 
\bigskip

Now let us consider $n+1$ particles $p_0,\dots,p_n$ at positions 
\[
p_{0,t}\leq p_{1,t}\leq\dots\leq p_{n,t}
\]
as functions of time $t$, initially all at the origin 0, and under the sorted dynamics where particles will always be indexed from left to right.   For $j=1,\dots,n$ we let $g_{j,t}=p_{j,t}-p_{j-1,t}$.  Thus $g_{j,0}=0$ for all $0$.  We also define $g_{0,t}\equiv\infty$, to be the gap ``to the left'' of the first particle. Assume that the process finishes within $T$ steps.

We will aim to show that it is unlikely that $p_n$ makes it to $x= Kn$ for some large constant $K$, by bounding the expected number of visits of $p_n$ to $x$.  
For particle $i$ and any $x\in \Z$, we define $\visits_i^{t_1}(x)$ to be the number of times $t\geq t_1$ that the particle $p_i$ \emph{arrives} at $x$ in time period $t\geq t_1$ of the process; this is the number of times $t\geq \max(1,t_1)$ that $p_{i,t-1}\neq x$ and $p_{i,t}=x$.  We write $\visits_i(x)=\visits_i^0(x)$.


For each particle $j$, we will consider its trajectory from its last coincidence with $p_{j-1}$ at some vertex $x_{j-1}$ until its next coincidence with $p_{j-1}$.  Observations \ref{bias} and \ref{ends} allow us to use Lemma \ref{walkvisits} to analyze the number of visits $p_j$ will make to a location $x$ before the next such coincidence.


In particular, as $g_{0,t}\equiv \infty$, we have for all $x\in \Z$ that the number $\visits_0(x)$ of times $t$ that the particle $p_0$ will ever land on $x$ (e.g., when $p_{0,t-1}\neq x$ and $p_{0,t}=x$) satisfies 
\begin{equation}
\label{0visitstail}
\E(\visits_0(x))\leq \rho^{\max(0,x)-c}
\end{equation}
by Lemma \ref{walkvisits}.
(Without using more about the process, we don't know whether the time $p_0$ actually spends at $x$ may greatly exceed this estimate; we are just bounding the number of steps when it arrives there.)

Now let $\xi_{i,t,x}$ be the indicator random variable
\[
\xi_{i,t,x}=\begin{cases}
    1&\text{ if }p_{i-1,t-1}=p_{i,t-1}=x\text{ and }p_{i,t}>p_{i-1,t}\\
    0& \text { otherwise,}
\end{cases}
\]
which is 1 for $i,t,x$ at times that $p_i$ has just moved to the right of $p_{i-1}$, after coinciding with it at $x$.  The plan of our proof is to use Lemma \ref{walkvisits} to analyze the trajectory of a particle $p_i$ each time it separates from $p_{i-1}$, until their next coincidence.  $\xi_{i,t,x}$ is the indicator random variable telling us that at time $t$, such a separation of $p_i$, $p_{i-1}$ has just happened at $x$.

Write $\fvisits_i^{t_1}(x)$ for the number of times $t_1\leq t\leq t'$ that $p_{i}$ arrives at $x$ ($p_{i,t-1}\neq x$, $p_{i,t}=x$), where $t'\leq \infty$ is minimum such that $p_{i,t'}=p_{i-1,t'}$.   This is precisely the quantity whose expectation we will control with Lemma \ref{walkvisits}.

Our goal now is to count arrivals of particle $j$ to a vertex $x$ using the indicator variables $\xi_{i,t,x}$ by considering the sum
\[
\sum_{x_{j-1}\in \Z}\sum_{t}\xi_{j,t,x_{j-1}} \fvisits_j^{t}(x).
\]
Given a value of $x$, this counts the number of times $t'$ such that:
\begin{itemize}
    \item $p_j$ arrives at $x$ at time $t'$,
    \item there is a $t<t'$ such that
    \begin{itemize}
        \item $p_j$ separated from $p_{j-1}$ at time $t'$, 
        \item there is no time $t''$ with $t<t''<t'$ when $p_j,p_{j'}$ coincide.
    \end{itemize}
\end{itemize}
Thus each arrival to $p_j$ is counted at most once in the sum, but the sum doesn't count arrivals where $p_j$ and $p_{j-1}$ move to $x$ together (either from $x-1$ or $x$).  This motivates the definition of $\svisits_i(x)$, which is the number of times $t$ when:
\begin{itemize}
\item $p_i$ arrives at $x$ (i.e., $p_{i,t-1}\neq x$, $p_{i,t}=x$) AND
\item at least one of $p_{i-1,t-1}=p_{i,t-1}$ and $p_{i-1,t}=p_{i,t}$ is false.
\end{itemize}
In particular, $\svisits_i(x)$ doesn't count the arrivals of $p_i$ to $x$ if it arrives there ``with'' $p_{i-1}$, and from the same vertex.
We have:
\begin{observation}
\[
\svisits_j(x)= \sum_{x_{j-1}\in \Z}\sum_{t}\xi_{j,t,x_{j-1}} \fvisits_j^{t}(x).
\]
\end{observation}
\begin{proof}
    Each contribution to $\svisits_j(x)$ comes from an arrival of particle $p_j$ to $x$ that occurs at a time $t'$ where we have that $p_{i-1,t'-1}\neq p_{i,t'-1}$ or else that $p_{i-1,t'-1}=p_{i,t'-1}$ and $p_{i-1,t'}\neq x$.  In the first case, there some maximum $t\leq t'-1$ such that $\xi_{j,t,y}=1$ for some $y$; in the second case, $\xi_{j,t',y}=1$ for $y=p_{i-1,t'-1}=p_{i,t'-1}$.  In either case, we can write $x_{j-1}=y$ for this $y$.   The visit at time $t'$ then contributes exactly 1 to the count $\fvisits_j^{t}(x)$; in particular, there can be no time $\tau\in [t,t')$ when $p_{i-1,\tau}=p_{i,\tau}$, because then either the visit at time $t'$ doesn't contribute to the count $\svisits_j(x)$, or else $t$ was not the largest $t\leq t'$ for which $\exists y,\, \xi_{j,t,y}=1$. 
\end{proof}

Now taking expectations, we have
\begin{multline}
\label{esvisits}
\E(\svisits_j(x))= \sum_{x_{j-1}\in \Z}\sum_{t}\E(\xi_{j,t,x_{j-1}} \fvisits_j^{t}(x))\\=
 \sum_{x_{j-1}\in \Z}\sum_{t}\Pr(\xi_{j,t,x_{j-1}}=1)\E\left(\fvisits_j^{t}(x)\mid \xi_{j,t,x_{j-1}}=1\right)
\end{multline}
Moreover, Lemma \ref{walkvisits} gives that:
\begin{observation}
\begin{equation}\label{condexp}
\E(\fvisits_j^{t}(x)\mid \xi_{j,t,x_{j-1}}=1)\leq \rho^{\max(x-x_{j-1},0)-c}.
\end{equation}
\end{observation}
\begin{proof}
    To handle conditioning on $\xi_{j,t,x_{j-1}=1}$, consider conditioning on any history of the process up to time $t$ compatible with $\xi_{j,t,x_{j-1}=1}$, e.g., any history up to time $t$ for which we had $p_{j-1,t-1}=p_{j,t-1}=x_{j-1}$ and $p_{j,t}=x_{j-1}+1,$ $p_{j-1,t}=x_{j-1}-1$.  Given any such history, we can couple the effect of the remaining process on the particle $p_j$ with the adversarial process that is the subject of Lemma \ref{walkvisits}, where ``ending'' the walk corresponds to the next time $\tau$ when $p_{j,\tau}=p_{j,\tau-1}$.   In particular, at any time $\tau>t$ when $p_j$ moves but $p_j$, $p_{j-1}$ have not yet coincided after time $t$, either:
    \begin{itemize}
        \item $p_{j-1,\tau}\leq p_{j,\tau}-2$, in which case there is at least a $\frac 3 4$ chance that $p_j$'s move is to the left, or else
        \item $p_{j-1,\tau}=p_{j,\tau}-1$, in which case there is at least probability $\frac 1 4$ that on the next step, we have $p_{j-1,\tau+1}=p_{j,\tau+1}$ (corresponding to ``ending'' the adversarial walk).
    \end{itemize}
    Lemma \ref{walkvisits} then gives the desired bound.  Note that at time $\tau=t$, the particle $p_j$ deterministically moves to $x_{j-1}+1$ since we are conditioning on $\xi_{j,t,x_{j-1}}=1$, so our coupled adversarial process actually begins with the token at 1, not 0 as stated in Lemma \ref{walkvisits}, and so we can use the choices for $\rho,c$ given in \eqref{n2}.
\end{proof}
Now \eqref{esvisits} and \eqref{condexp} give that
\begin{multline}
    \label{esex}
    \E(\svisits_j(x))\leq \sum_{x_{j-1}\in \Z}\rho^{\max(x-x_{j-1},0)-c}\sum_{t}\Pr(\xi_{j,t,x_{j-1}}=1)\\=
    \sum_{x_{j-1}\in \Z}\E\left(\sum_{t}\xi_{j,t,x_{j-1}}\right)\rho^{\max(x-x_{j-1},0)-c}.
\end{multline}

Now we will use:
\begin{observation}
\begin{equation}
\label{sumxi}\sum_{t}\xi_{j,t,x} \leq\1_{\substack{x=0}}+\visits_{j-1}(x)\leq \1_{\substack{x=0}}+\sum_{i< j}\svisits_i(x).
\end{equation}
\end{observation}
\begin{proof}
The first inequality follows from the fact that each time $t$ that $p_j$ and $p_{j-1}$ coincide at $x$ at time $t-1$ and then are separated, $p_{j-1}$  leaves $x$ (in particular, $p_{j-1,t}=x-1$, $p_{j,t}=x+1$), and so we can associate each such separation to the last time $p_{j-1}$ arrived at $x$, except for the very first separation from $x=0$, which is accounted for by the $\1_{x=0}$ term.  The second inequality is because for each arrival of $p_{j-1}$ at $x$ at time $t>0$, there is some smallest $i<j$ such that $p_{i}$ is arriving at $x$ at time $t$, and the arrival is counted by $\svisits_i(x)$.\\




\end{proof}

Applying \eqref{sumxi} to \eqref{esex}, we have that
\begin{multline}\label{soloarrivals}
\E(\svisits_j(x))\leq \sum_{x_{j-1}\in \Z}\E\left(\sum_{t}\xi_{j,t,x_{j-1}}\right)\rho^{\max(0,x-x_{j-1})-c}\\\leq
 \sum_{x_{j-1}\in \Z}\left(\1_{x_{j-1}=0}+\sum_{i<j}\E(\svisits_i(x_{j-1}))\right
)\rho^{\max(0,x-x_{j-1})-c}
\end{multline}

Now let $X_j(x)$ denote the set of all sequences $0=x_{-1},x_0,x_1,x_2,\dots,x_j=x$, for $x_i\in \Z$.
For some sequence $\sigma\in X_j(x)$, we let 
\[
\phi(\sigma)=\sum_{i=0}^{j} \max(x_{i-1}-x_{i},0).
\]
Using \eqref{soloarrivals}, we aim to prove inductively that 
\begin{equation}
    \label{toprove}
    \E(\svisits_j(x))\leq (j+1)\sum_{\sigma\in X_{j}(x)}\rho^{x+\phi(\sigma)-jb},
\end{equation}
for some constant $b\geq c$.  Equation \eqref{0visitstail} serves as the base case for $j=0$.  For the sake of intuition, note that the $x_0,x_1,\dots,x_{j-1}$ corresponds to a sequence for which there exist times $t_0\leq t_1\leq\dots\leq t_{j-1}$ such that for each $1\leq i\leq j-2$, $p_{i-1}$ and $p_i$ coincided at $x_{i-1}$ at time $t_{i-1}$, and did not coincide again before time $t_{i}$.  Sequences $\sigma$ with large $\phi(\sigma)$ are wasteful, in the sense that progress to the right made by one particle is lost by the next one, which ends up limiting the contribution of such sequences to $\svisits_j(x)$, as will be captured by \eqref{toprove}.

Note that the induction hypothesis \eqref{toprove} applied for $i<j$
gives that
\begin{multline}
\sum_{i<j}\E(\svisits_i(x))\leq \sum_{i<j}(i+1)\sum_{\sigma\in X_{i}(x)}\rho^{x+\phi(\sigma)-ib}\\
=\sum_{i<j}(i+1)\rho^{x-ib}\sum_{\sigma\in X_i(x)}\rho^{\phi(\sigma)}
\leq
j\sum_{\sigma\in X_{j-1}(x)} \rho^{x+\phi(\sigma)-(j-1)b}\sum_{k=0}^{j-1}\rho^{kb}\\
\leq
\frac{j}{1-\rho^b}\sum_{\sigma\in X_{j-1}(x)} \rho^{x+\phi(\sigma)-(j-1)b}.
\end{multline}

In particular, returning to \eqref{soloarrivals}, we get by induction and
 and the fact that 
\[
\sum_{\sigma\in  X_i(x)} \rho^{\phi(\sigma)}\geq \sum_{\sigma\in X_{i-1}(x)} \rho^{\phi(\sigma)}\quad\text{for all $i$.}
\]
we get that
\begin{multline}
\E(\svisits_j(x))\\\leq 
 \sum_{x_{j-1}\in \Z}\left(\1_{x_{j-1}=0}+\sum_{i<j}\E(\svisits_i(x_{j-1}))\right)\rho^{\max(0,x-x_{j-1})-c}\\\leq
 \sum_{x_{j-1}\in \Z}  \left(\1_{x_{j-1}=0}+\frac{j}{1-\rho^b}\sum_{\sigma\in X_{j-1}(x)}\rho^{x_{j-1}+\phi(\sigma)-(j-1)b}\right) \rho^{\max(0,x-x_{j-1})-c}\\
 =\sum_{x_{j-1}\in \Z}  \left(\1_{x_{j-1}=0}\rho^{x-c}+\frac{j}{1-\rho^b}\sum_{\sigma\in X_{j}(x)}\rho^{x+\phi(\sigma)-(j-1)b-c}\right)\\
 \leq \1_{x_{j-1}=0}\rho^{x-c}+\frac{j}{1-\rho^b}\sum_{\sigma\in X_{j}(x)}\rho^{x+\phi(\sigma)-(j-1)b-c}\\
 \leq (j+1)\sum_{\sigma\in X_{j}(x)}\rho^{x+\phi(\sigma)-jb},
\end{multline}
by choosing $b=c+1$ as in \eqref{n1}, \eqref{n2}, giving
\[
\rho^{c-b}\geq \frac{1}{1-\rho^b}.
\]
This completes the inductive proof of \eqref{toprove}.

\bigskip

In particular, the expected number of visits of $p_n$ to $Kn$ is at most
\begin{multline}
\label{endgame}
\sum_{j\leq n} \E(\svisits_j(Kn))\leq \sum_{j\leq n}(j+1)\sum_{\sigma\in {X_{j}(Kn)}}\rho^{Kn+\phi(\sigma)-jb}\\\leq 
n^2 \rho^{(K-b)n} \sum_{\sigma\in {X_{n}(Kn)}}\rho^{\phi(\sigma)},
\end{multline}
and so it remains to bound the sum 
\[
\sum_{\sigma\in {X_{n}(Kn)}}\rho^{\phi(\sigma)}.
\]
For the sake of intuition, note first that the contribution to the sum from terms for which $\phi(\sigma)=0$ is precisely
\[
\binom{Kn+n-1}{n-1},
\]
whose product with $\rho^{Kn}$ goes to 0 very quickly.  In particular, by Markov's inequality and \eqref{endgame} the proof of Theorem \ref{disperse} is complete if we can show that 
\begin{equation}\label{finish}
n^2\rho^{(K-b)n}\sum_{\sigma\in X_{n}(Kn)}\rho^{\phi(\sigma)}=o(1).
\end{equation}

To this end, let $m=Kn$ and
\[
X_n^f(m)=\{x\in X_n(m)\mid \phi(x)=f\},\quad A_f=|X_n^f(m)|.
\]
Note that 
\begin{equation}\label{A0}
A_0=\binom{m+n}{n}\leq (e(K+1))^n.
\end{equation}

We have
\[
A_f=\sum_{\ell\geq 0}\sum_{|S|=\ell} \sum_{\substack{\sum_{i\in S} f_i=f\\\forall i,\,f_i\geq 1}} |\{x:x_{i}=x_{i-1}+g_i\}|
\]
where $g_i=-f_i,i\in S$ and $g_i\geq 0$ for $i\notin S$ and
\[
\sum_{i\notin S}g_i=m+f.
\]
So for $f\geq 1$ we have 
\begin{multline*}
A_f\leq \sum_{\ell=1}^{\min\{f,n\}}\sum_{|S|=\ell} \binom{f-1}{\ell-1} \binom{m+f+n-\ell}{n-\ell}=\\
\sum_{\ell=1}^{\min\{f,n\}}\binom{n}{\ell}\binom{f-1}{\ell-1} \binom{m+f+n-\ell}{n-\ell}\leq
2^n \sum_{\ell=1}^{\min\{f,n\}}\binom{f-1}{\ell-1} \binom{m+f+n-\ell}{n-\ell}.
\end{multline*}

Now consider cases as follows:\\
\noindent \textbf{Case 1:} $f\leq Kn$.\\
Then as $m=Kn$ we have
\[
A_f\leq 2^n n \binom{Kn}{n}\binom{m+Kn+n}{n}\leq 2^n n  (eK)^n (e(2K+1))^n.
\]
\textbf{Case 2:} $f> Kn$.\\
Then we have
\[
A_f\leq 2^n n \max_{\ell\leq n}\left(\binom{f}{\ell}\binom{m+f+n}{n}\right)\leq
2^n n \left(\frac{fe}{n}\right)^n\left(\frac{2ef}{n}\right)^n
\]
In particular we get that
\begin{multline}\label{workingit}
\sum_{\sigma\in X_n(m)}\rho^{\phi(\sigma)}=\sum_{f\geq 0}\sum_{\sigma\in X_{n}^f(m)}\rho^f
\leq
\sum_{f=0}^{Kn}|A_f|+
\sum_{f=Kn}^{n^2} |A_f|\rho^f\\\leq
(Kn+1)(10K)^{4n}+
\frac{20^n}{n^{2n}}\sum_{f=Kn}^{n^2}\rho^f f^{2n}
\end{multline}
In the last sum, the ratio of the $(f+1)$st term to the $f$th term is
\[
\rho\left(1+\frac{1}{f}\right)^{2n}\leq \rho e^{2n/f},
\]
which is at most 1 for 
\[
f\geq \frac{2n}{\ln(\rho^{-1})}.
\]
For example, for $\rho=\sqrt{3}/{2}$, if $K\geq 15$, so that $f\geq 15n$ in the last sum in \eqref{workingit}, we have
\[
\sum_{\sigma\in X_n(m)}\rho^{\phi(\sigma)}\leq 
(Kn+1)(10K)^{4n}+
\frac{20^nn^2}{n^{2n}}\rho^{Kn} (Kn)^{2n}\leq K^{10n}
\]
for large $n$, and this proves \eqref{finish} for large $n$ and $K>b$ (e.g., we can take $K=62$).\qed

\section{Remaining proofs}
\label{remaining}
It remains to prove Observation \ref{ends} and Lemma \ref{walkvisits}.

\begin{proof}[Proof of Observation \ref{ends}]
Let $x=p_{i-1,t-1},a=N_{x,t-1}$ and $b=N_{x+1,t-1}$.
  If $x+1$ topples at time $t-1$ but $x$ does not, we have $\Pr(p_{i,t}=p_{i,t-1})=(1-2^{-b})$.  
  
  If both topple, let $k=\mathrm{Binom}(a,\tfrac 1 2),\ell=\mathrm{Binom}(b,\tfrac 1 2)$ be the random variables counting the the number of particles that move left from the pile at $x$ and $x+1$, respectively.  Let $L_{x-1}$, $L_x$ and $L_{x+1}$ denote the number of particles at vertices $\leq x-1$, $\leq x$, and $\leq x+1$ after these two topples.  We have
  \[
  L_{x-1}=i-1-(a-k),\quad L_x=i-1+\ell-(a-k),\quad L_{x+1}\geq i-1+\ell.
  \]
  So if $\ell>a-k$ and $a-k>0$, then $L_{x-1}<i-1$ and $L_x\geq i$, so that $p_i$ and $p_{i-1}$ coincide at $x$.  On the other hand, if $a-k>\ell$ and $\ell\geq 1$, then $L_x\leq i-2$ and $L_{x+1}\geq i$, so that $p_i$ and $p_{i-1}$ coincide at $x+1$.  

  So $p_i$ and $p_{i-1}$ will coincide if $k<a$ and $\ell>0$ and $a\neq k+\ell$.  For $a=b=2$ one can check that this has probability exactly $\frac 1 4$.  More generally, we have
    \begin{multline}
    \Pr(k<a)\Pr(\ell>0)\Pr(k+\ell\neq a\mid k<a,\ell>0)\\=
    (1-2^{-a})(1-2^{-b})\left(1-\sum_{k=0}^{a-1}\frac{\binom{a}{k}}{2^a-1}\frac{\binom{b}{a-k}}{2^b-1}\right)\\=
   (1-2^{-a})(1-2^{-b})\left(1-\frac{\binom{a+b}{a}-1}{(2^a-1)(2^b-1)}\right)= (1-2^{-a})(1-2^{-b})-\frac{\binom{a+b}{a}-1}{2^{a+b}}.
\end{multline}
For $a+b\geq 4$ we have that $\left(\binom{a+b}{a}-1\right)/2^{a+b}\leq \frac{5}{16}$ which gives that
\begin{multline}
    \Pr(k<a)\Pr(\ell>0)\Pr(k+\ell\neq a\mid k<a,\ell>0)\\ \geq
    (1-2^{-a})(1-2^{-b})-\frac{5}{16},
\end{multline}
which is then minimized by minimizing $a,b$.  Using $a=b=2$ gives the claimed bound of $\frac 1 4$.
\end{proof}

\begin{proof}[Proof of Lemma \ref{walkvisits}]
First note that suffices to prove the the lemma for $x\geq 0$, since for any fixed $x<0$, we can simply start the analysis once the process has reached $x$.  So we will assume $x\geq 0$; in particular, we can take $\max(0,x)=x$ in \eqref{adversary}.

Ahead of time, we fix 4 separate random streams that will be used for:
\begin{enumerate}[(A)]
\item \label{stepstream} Steps when the adversary chooses option \ref{step},
\item \label{epstream} The choice of $\epsilon_t$, when the adversary chooses option \ref{term},
\item \label{termstream} Any other randomness used by the adversary when choosing option \ref{term}, and
\item \label{choicestream} Any randomness used by the adversary to determine which option to choose on a given step.
\end{enumerate}
Each of these stream can be implemented as a sequence of independent uniform $U[0,1]$ random variables; when making a choice that requires randomness, the adversary observes the next unobserved uniform in the appropriate sequence, and can use the value of that random variable as input to the decision.  In particular, preselecting the random stream for \ref{epstream} is equivalent to choosing a geometric random variable $\bar\kappa_2\geq 1$ of parameter $r$ so that the adversary will observe $\epsilon_t=1$ if and only if they are choosing option \ref{term} for the $\bar\kappa_2$nd time.  

Now fix some strategy $\sigma$ of the adversary and let $S_0,S_1,\dots$ denote the (infinite) random walk the strategy would produce under an alternative scenario where the random stream for \ref{epstream} is not used, and instead on every step $\sigma$ chooses option \ref{term}, the strategy observes $\epsilon_t=0$.  Here $S_k=\sum_{\ell\leq k}\xi_k$, where $\xi_k\in \{-1,0,+1\}$ is the change in position at the $k$th step of the walk.  Note that by our pre-selection of randomness, the walk produced by $\sigma$ when using all the pre-selected randomness (i.e., not always just observing $\epsilon_t=0$) will be either this whole walk or a finite initial segment of it, depending on the adversary's choices. 

    But considering still what $\sigma$ produces with the the pre-selected randomness for \ref{stepstream}, \ref{termstream}, \ref{choicestream} but always observing $\epsilon_t=0$, when choosing option \ref{term}, define now $\xi_{i_1},\xi_{i_2},\dots$ to be the subsequence of the $\xi_i$ corresponding to steps when the adversary chooses option \ref{step}.  Note that each for each $j$, we have that
    \[
    \Pr\left(\xi_{i_j}=-1\mid \xi_{i_{j-1}},\xi_{i_{j-2}},\dots,\xi_{i_1}\right)\geq p.
    \]
    Thus $T_k=\sum_{j\leq k} \xi_{i_j}$ defines a walk $T_0,T_1,\dots$ on $\Z$ that at each step moves to the left with probability at least $p.$ In particular, we have that there is a $\lambda<1$ such that
    \begin{equation}\label{forT}
    \forall x\in \Z, \E\Big(\big|\{n\mid T_n=x\}\big|\Big)\leq \lambda^{x-d}.
    \end{equation}
    Indeed,
    \begin{align*}
\E\Big(\big|\{n\mid T_n=x\}\big|\Big)&=\sum_{n\geq x}\binom{n}{(n+x)/2}p^{(n-x)/2}(1-p)^{(n+x)/2}  \\
&=\left(\frac{1-p}{p}\right)^{x/2}\sum_{n\geq x}\binom{n}{(n+x)/2}(p(1-p))^{n/2}\\
&\leq \left(\frac{1-p}{p}\right)^{x/2}\sum_{n\geq x}(4p(1-p))^{n/2}.\\
    \end{align*}
In particular, we can take any
\[
\lambda\leq \sqrt{(1-p)/p)},\quad d\geq -\log_p\left(\frac{1}{1-\sqrt{4p(1-p)}}\right).
\]
So, for $p=3/4$, we can take $\lambda=1/\sqrt{3}$, $d=7$. Now let $\kappa_2\geq 0$ denote the number of times the strategy $\sigma$ chooses option \ref{term} when using all the pre-selected randomness (i.e., not when just observing $\epsilon_t=0$ deterministically); we have $\kappa_2\leq \bar\kappa_2$ and so $\kappa_2<\infty$ with probability 1.  (Note that the adversary is not required to use option 2 $\bar\kappa$ times.)  The actual walk produced by $\sigma$ using all the preselected sources of randomness is either the whole walk $S_n$ if $\kappa_2<\bar\kappa_2$, and an initial segment of it if $\kappa_2=\bar\kappa_2$ (ending at the step where the adversary chooses option \ref{term} for the $\bar\kappa_2$nd time).  

      Observe that (very crudely), we have that 
      \begin{equation}
    \big|\{n\mid S_n=x\}\big|\leq \sum_{y=x-\kappa_2}^{x+\kappa_2} \kappa_2\big|\{n\mid T_n=y\}\big|.
      \end{equation}
      Indeed, if we define $0\leq \nu(n)\leq n$ to be the number of times $\leq n$ the strategy $\sigma$ chooses option 1, then we have that $|S_{n}-T_{\nu(n)}|\leq \kappa_2$;    Moreover, $n\mapsto \nu(n)$ is at most a $\kappa_2$-to-1 mapping, so that any visit of $S_n$ to $x$ is one of $\kappa_2$ visits of $\{S_{n}\}$ to $x$ that correspond to a visit of $T_{\nu(n)}$ to a vertex within $\kappa_2$ of $x$.

      We thus have the inequality
    \begin{multline}
    \big|\{n\mid S_n=x\}\big|\leq \sum_{y=x-\kappa_2}^{x+\kappa_2} \kappa_2\big|\{n\mid T_n=y\}\big|
    =\sum_{y\in \Z} \kappa_2 \1_{\kappa_2\geq |y-x|} \big|\{n\mid T_n=y\}\big|
    \\\leq \sum_{y\in \Z}\bar \kappa_2 \1_{\bar \kappa_2\geq |y-x|} \big|\{n\mid T_n=y\}\big|.
    \end{multline}
    Now, by our construction, the infinite walks $\{S_n\}$ and $\{T_n\}$ are independent of $\bar\kappa_2$.  In particular, together with \eqref{forT}, and using  $\rho\geq\max(\sqrt{1-r},\lambda)$, we obtain
    \begin{multline}
    \E\left(\big|\{n\mid S_n=x\}\big|\right)\leq \sum_{y\in \Z}\E\Big(\bar\kappa_2\1_{\bar \kappa_2\geq |y-x|} \big|\{n\mid T_n=y\}\big|\Big)\\=
    \sum_{y\in \Z}\E\big(\bar\kappa_2\big| \bar\kappa_2\geq |y-x|\big)\Pr\big(\bar\kappa_2\geq |y-x|\big)\E\Big(\big|\{n\mid T_n=y\}\big|\Big)
    \\\leq
    \sum_{y\in \Z}\left(|y-x|-1+\frac 1 r\right)(1-r)^{|y-x|-1}
    \lambda^{y-d}
    \\
    =\sum_{k=-\infty}^\infty \left(|k|-1+\frac 1 r\right) (1-r)^{|k|-1}\lambda^{x+k-d}
    \\
    \leq
    \sum_{k=-\infty}^\infty \left(|k|-1+\frac 1 r\right) \rho^{2|k|-2}\rho^{x+k-d}
    \\\leq
    \sum_{k=-\infty}^\infty \left(|k|-1+\frac 1 r\right) \rho^{x+|k|-d-2}
    \\\leq
    2\rho^{x-d-2} \sum_{k=0}^\infty \left(k-1+\frac 1 r\right) \rho^{k}
    \\
    =2\rho^{x-d-2}\left(\frac{2r\rho-r-\rho+1}{r(1-\rho)^2}\right)
    \\
    \leq \frac{3}{r(1-\rho)^2}\rho^{x-d-2}\leq \rho^{x-c}
    \end{multline}
    for a constant $c$.  In particular, for the special case where $r=\frac 1 4$, $p=\frac 3 4$, we can take
    \[
    \lambda=\frac 1 {\sqrt 3}\quad d=7\quad \rho=\frac{\sqrt{3}}{2}\quad c=59.
    \]
\end{proof}

\end{document}